\documentclass{amsart}
\usepackage[all]{xy}
\usepackage{verbatim}
\usepackage{color}
\usepackage{amsthm}
\usepackage{amssymb}
\usepackage[colorlinks=true]{hyperref}
\usepackage{footmisc}
\numberwithin{equation}{section}

\newtheorem{theorem}[equation]{Theorem}
\newtheorem*{theorem*}{Theorem}

\newtheorem*{conjecture*}{Mamma Conjecture}
\newtheorem*{conjecture1*}{Mamma Conjecture (revisited)}

\newtheorem*{corollary*}{Corollary}

\theoremstyle{remark}

\newtheorem{example}[equation]{Example}

\theoremstyle{remark}
\newtheorem{remark}[equation]{Remark}

\newcommand{\cA}{{\mathcal A}}
\newcommand{\cB}{{\mathcal B}}
\newcommand{\cC}{{\mathcal C}}

\newcommand{\cF}{{\mathcal F}}

\newcommand{\cN}{{\mathcal N}}
\newcommand{\cO}{{\mathcal O}}

\newcommand{\cW}{{\mathcal W}}
\newcommand{\cX}{{\mathcal X}}
\newcommand{\cY}{{\mathcal Y}}
\newcommand{\cZ}{{\mathcal Z}}

\newcommand{\bbA}{\mathbb{A}}
\newcommand{\bbB}{\mathbb{B}}
\newcommand{\bbC}{\mathbb{C}}

\newcommand{\bbP}{\mathbb{P}}

\newcommand{\bbQ}{\mathbb{Q}}
\newcommand{\bbZ}{\mathbb{Z}}

\DeclareMathOperator{\id}{id}

\DeclareMathOperator{\NChow}{NChow} % category of noncommutative Chow motives
\DeclareMathOperator{\NHom}{NHom} % category of noncommutative Chow motives
\DeclareMathOperator{\NNum}{NNum} % category of noncommutative numerical motives

\newcommand{\dgcat}{\mathrm{dgcat}} % codimension 
\newcommand{\perf}{\mathrm{perf}}

\newcommand{\Chow}{\mathrm{Chow}}

\newcommand{\dg}{\mathrm{dg}}

\newcommand{\Hom}{\mathrm{Hom}}
\newcommand{\End}{\mathrm{End}}

\newcommand{\op}{\mathrm{op}}

\newcommand{\too}{\longrightarrow}

\newcommand{\ie}{\textsl{i.e.}\ }

\let\oldmarginpar\marginpar
\def\marginpar#1{\oldmarginpar{\tiny #1}}

\begin{document}

\title[A note on Grothendieck's standard conjectures of types $C^+$ and $D$]{A note on Grothendieck's \\standard conjectures of type $C^+$ and $D$}
\author{Gon{\c c}alo~Tabuada}

\address{Gon{\c c}alo Tabuada, Department of Mathematics, MIT, Cambridge, MA 02139, USA}
\email{tabuada@math.mit.edu}
\urladdr{http://math.mit.edu/~tabuada}
\thanks{The author was partially supported by a NSF CAREER Award}

\subjclass[2010]{14A22, 14C15, 14M12, 18D20, 18E30}
\date{\today}

%\keywords{\Goncalo{[At the end]}}
%
\abstract{Grothendieck conjectured in the sixties that the even K\"unneth projector (with respect to a Weil cohomology theory) is algebraic and that the homological equivalence relation on algebraic cycles coincides with the numerical equivalence relation. In this note we extend these celebrated conjectures from smooth projective schemes to the broad setting of smooth proper dg categories. As an application, we prove that Grothendieck's original conjectures are invariant under homological projective duality. This leads to a proof of Grothendieck's conjectures in the case of intersections of quadrics, linear sections of determinantal varieties, and intersections of bilinear divisors. Along the way, we prove also the case of quadric fibrations.}}

\maketitle
\vskip-\baselineskip
\vskip-\baselineskip
\vskip-\baselineskip
%\tableofcontents

%\bigskip

%\medskip
%-------------------------------------------------------------------------------
\section{Introduction and statement of results}
%-------------------------------------------------------------------------------
Let $k$ be a base field of characteristic zero. Given a smooth projective $k$-scheme $X$ and a Weil cohomology theory $H^\ast$, let us denote by $\pi_X^i\colon H^\ast(X) \to H^\ast(X)$ the $i^{\mathrm{th}}$ K\"unneth projector, by $Z^\ast(X)_\bbQ$ the $\bbQ$-vector space of algebraic cycles on $X$, and by $Z^\ast(X)_\bbQ/_{\!\sim \mathrm{hom}}$ and $Z^\ast(X)_\bbQ/_{\!\sim \mathrm{num}}$ the quotients with respect to the homological and numerical equivalence relations, respectively. Following Grothendieck \cite{Grothendieck} (see also Kleiman \cite{Kleim1, Kleim}), the standard conjecture\footnote{The standard conjecture of type $C^+$ is also usually known as the {\em sign conjecture}. If the even K\"unneth projector is algebraic, then the odd K\"unneth projector $\pi^-_X:=\sum_i \pi_X^{2i+1}$ is also algebraic.} of type $C^+$, denote by $C^+(X)$, asserts that the even K\"unneth projector $\pi^+_X:=\sum_i \pi^{2i}_X$ is algebraic, and the standard conjecture of type $D$, denoted by $D(X)$, asserts that $Z^\ast(X)_\bbQ/_{\!\sim \mathrm{hom}}=Z^\ast(X)_\bbQ/_{\!\sim \mathrm{num}}$. Thanks to the work of Kleiman \cite{Kleim} and Lieberman \cite{Lieberman}, and to the fact that $D(X\times X) \Rightarrow C^+(X)$ (see \cite[Thm.~5.4.2.1]{Andre}), the conjecture $C^+(X)$, resp. $D(X)$, holds in the case where $X$ is of dimension $\leq 2$, resp. $\leq 4$, and also for abelian varieties. In addition to these cases, the aforementioned important conjectures remain~wide~open.

A {\em dg category} $\cA$ is a category enriched over complexes of $k$-vector spaces; see \S\ref{sub:dg}. Every (dg) $k$-algebra $A$ gives naturally rise to a dg category with a single object. Another source of examples is provided by schemes since the category of perfect complexes $\perf(X)$ of every quasi-compact quasi-separated $k$-scheme $X$ admits a canonical dg enhancement\footnote{When $X$ is quasi-projective this dg enhancement is unique; see Lunts-Orlov \cite[Thm. 2.12]{LO}.} $\perf_\dg(X)$. As explained in \S\ref{sub:C}-\ref{sub:D}, given a smooth proper dg category $\cA$ in the sense of Kontsevich, the standard conjectures of type $C^+$ and $D$ admit noncommutative analogues $C^+_{\mathrm{nc}}(\cA)$ and $D_{\mathrm{nc}}(\cA)$, respectively.
% let us denote by $K_0(\cA)_\bbQ$ the $\bbQ$-linearization of the Grothendieck group of $\cA$. As explained in \S\ref{sub:numerical}, this $\bbQ$-vector space comes equipped with an homological equivalence relation ${}_{\!\sim\mathrm{hom}}$ and with a numerical equivalence relation ${}_{\!\sim\mathrm{num}}$. Following \cite{JEMS}, the noncommutative standard conjecture of type $D$, denoted by $D_{\mathrm{nc}}(\cA)$, asserts that $K_0(\cA)_\bbQ/_{\!\sim \mathrm{hom}}=K_0(\cA)_\bbQ/_{\!\sim \mathrm{num}}$.
\begin{theorem}\label{thm:main}
Given a smooth projective $k$-scheme $X$, we have the following equivalences\footnote{Thanks to Theorem \ref{thm:main}, the assumption that Grothendieck's standard conjectures of type $C^+$ and $D$ hold for every smooth projective $k$-scheme can now be removed from \cite[Thm.~1.7]{JEMS}.} of conjectures $C^+(X)\Leftrightarrow C^+_{\mathrm{nc}}(\perf_\dg(X))$ and $D(X)\Leftrightarrow D_{\mathrm{nc}}(\perf_\dg(X))$.
\end{theorem}
Theorem \ref{thm:main} extends Grothendieck's standard conjectures of type $C^+$ and $D$ from schemes to dg categories. Making use of this noncommutative viewpoint, we now prove Grothendieck's original conjectures in the case of quadric fibrations:
\begin{theorem}[Quadric fibrations]\label{thm:app1}
Let $q\colon Q\to S$ be a flat quadric fibration of relative dimension $d$. Whenever the dimension of $S$ is $\leq 2$, resp. $\leq 4$, $d$ is even,~and the discriminant divisor of $q$ is smooth, the conjecture $C^+(Q)$, resp. $D(Q)$, holds.
\end{theorem}
\begin{remark}
A ``geometric'' proof of Theorem \ref{thm:app1} can be obtained by combining the aforementioned work of Kleiman and Lieberman, with Vial's computation (see \cite[Thm.~4.2 and Cor.~4.4]{Vial}) of the rational Chow motive of $Q$. 
\end{remark}
Making use of Theorem \ref{thm:main}, we now prove that Grothendieck's conjectures are invariant under homological projective duality (=HPD). Let $X$ be a smooth projective $k$-scheme equipped with a line bundle $\cO_X(1)$; we write $X \to \bbP(V)$ for the associated morphism where $V:=H^0(X,\cO_X(1))^\ast$. Assume that the triangulated category $\perf(X)$ admits a Lefschetz decomposition $\langle \bbA_0, \bbA_1(1), \ldots, \bbA_{i-1}(i-1)\rangle$ with respect to $\cO_X(1)$ in the sense of Kuznetsov \cite[Def.~4.1]{KuznetsovHPD}. Following \cite[Def.~6.1]{KuznetsovHPD}, let $Y$ be the HP-dual of $X$, $\cO_Y(1)$ the HP-dual line bundle, and $Y\to \bbP(V^\ast)$ the morphism associated to $\cO_Y(1)$. Given a generic linear subspace $L \subset V^\ast$, consider the linear sections $X_L:=X\times_{\bbP(V^\ast)} \bbP(L^\perp)$ and $Y_L:=Y \times_{\bbP(V)} \bbP(L)$. 
\begin{theorem}[HPD-invariance]\label{thm:HPD}
Let $X$ and $Y$ be as above. Assume that $X_L$ and $Y_L$ are smooth, that $\mathrm{dim}(X_L)=\mathrm{dim}(X) -\mathrm{dim}(L)$, that $\mathrm{dim}(Y_L)=\mathrm{dim}(Y)- \mathrm{dim}(L^\perp)$, and that the conjecture $C^+_{\mathrm{nc}}(\bbA_0^{\dg})$, resp. $D_{\mathrm{nc}}(\bbA_0^{\dg})$, holds, where $\bbA_0^\dg$ stands for the dg enhancement of $\bbA_0$ induced from $\perf_\dg(X)$. Under these assumptions, we have the equivalence  $C^+(X_L)\Leftrightarrow C^+(Y_L)$, resp. $D(X_L)\Leftrightarrow D(Y_L)$. % Consequently, the latter conjecture holds when $X_L$ or $Y_L$ is of dimension $\leq 4$.
\end{theorem}
\begin{remark}
\begin{itemize}
\item[(i)] Conjectures $C^+_{\mathrm{nc}}(\bbA_0^{\dg})$ and $D_{\mathrm{nc}}(\bbA_0^{\dg})$ hold, in particular, whenever the triangulated category $\bbA_0$ admits a full exceptional collection.
\item[(ii)] Theorem \ref{thm:HPD} holds more generally when $X$ (or $Y$) is singular. In this case, we need to replace $X$ by a noncommutative resolution of singularities $\perf_\dg(X;\cF)$ and conjecture $C^+(X_L)$, resp. $D(X_L)$, by conjecture $C^+_{\mathrm{nc}}(\perf_\dg(X_L;\cF_L))$, resp. $D_{\mathrm{nc}}(\perf_\dg(X_L;\cF_L))$, where $\cF_L$ stands for the restriction of $\cF$ to $X_L$; consult \cite[\S2.4]{ICM-Kuznetsov} for further details. 
\end{itemize}
\end{remark}
To the best of the author's knowledge, Theorem \ref{thm:HPD} is new in the literature. As a first application, it provides us with an alternative (noncommutative) formulation of Grothendieck's original conjectures. Here are two ``antipodal'' examples; many more can be found at the survey \cite{ICM-Kuznetsov}.
\begin{example}[Veronese-Clifford duality]
Let $W$ be a $k$-vector space of dimension $d$, and $X$ the associated projective space $\bbP(W)$ equipped with the double Veronese embedding $\bbP(W) \to \bbP(S^2W)$. By construction, we have a flat quadric fibration $q\colon Q \to \bbP(S^2W^\ast)$, where $Q$ stands for the universal quadric in $\bbP(W)$. As proved in \cite[Thm.~5.4]{Kuznetsov-quadrics}, the HP-dual $Y$ of $X$ is given by $\perf_\dg(\bbP(S^2W^\ast);\cF)$, where $\cF$ stands for the sheaf of even Clifford algebras associated to $q$. Moreover, given a generic linear subspace $L \subset S^2W^\ast$, the linear section $X_L$ corresponds to the (smooth) intersection of the $\mathrm{dim}(L)$ quadric hypersurfaces in $\bbP(W)$ parametrized by $L$, and $Y_L$ is given by $\perf_\dg(\bbP(L);\cF_L)$. Making use of Theorem \ref{thm:HPD}, we hence conclude that $C^+(X_L) \Leftrightarrow C^+_{\mathrm{nc}}(\perf_\dg(\bbP(L);\cF_L))$ and $D(X_L) \Leftrightarrow D_{\mathrm{nc}}(\perf_\dg(\bbP(L);\cF_L))$.
\end{example}
By solving the preceding noncommutative standard conjectures, we hence prove Grothendieck's original standard conjectures in the case of intersections of quadrics:
\begin{theorem}[Intersections of quadrics]\label{thm:app2}
Whenever the dimension of $L$ is $\leq 3$, resp. $ \leq 5$, $d$ is even, and the discriminant division of $q$ is smooth, the conjecture $C^+_{\mathrm{nc}}(\perf_\dg(\bbP(L);\cF_L))$, resp. $D_{\mathrm{nc}}(\perf_\dg(\bbP(L);\cF_L))$, holds. Consequently, Grothendieck's original standard conjecture $C^+(X_L)$, resp. $D(X_L)$, also holds. 
\end{theorem}
\begin{remark}
As mentioned by Grothendieck at \cite[page~197]{Grothendieck}, the standard conjecture of Lefschetz type $B(X)$ holds for smooth complete intersections. Since this conjecture implies the standard conjectures of type $C^+$ and $D$ (the implication $B(X) \Rightarrow D(X)$ uses in an essential way the Hodge index theorem; see \cite[Thm.~4.1 and Prop.~5.1]{Kleim1}), we hence obtain an alternative ``geometric'' proof of Theorem \ref{thm:app2}.
\end{remark}
\begin{example}[Grassmannian-Pfaffian duality]
Let $W$ be a $k$-vector space of dimension $6$, and $X$ the associated Grassmannian $\mathrm{Gr}(2,W)$ equipped with the Pl\"ucker embedding $\mathrm{Gr}(2,W) \to \bbP(\bigwedge^2 W)$. As proved in \cite[Thm.~1]{Lines}, the HP-dual $Y$ of $X$ is given by $\perf_\dg(\mathrm{Pf}(4,W^\ast);\cF)$, where $\mathrm{Pf}(4,W^\ast) \subset \bbP(\bigwedge^2W^\ast)$ is the singular Pfaffian variety and $\cF$ a certain coherent sheaf of algebras\footnote{The sheaf of algebras $\cF$ is isomorphic to a matrix algebra on the smooth locus of $\mathrm{Pf}(4,W^\ast)$. Consequently, whenever $\mathrm{Pf}(4,W^\ast)_L$ is contained in the smooth locus of $\mathrm{Pf}(4,W^\ast)$, we conclude that the conjectures $C^+_{\mathrm{nc}}(\perf_\dg(\mathrm{Pf}(4,W^\ast)_L; \cF_L))$ and $D_{\mathrm{nc}}(\perf_\dg(\mathrm{Pf}(4,W^\ast)_L; \cF_L))$ are equivalent to $C^+(\mathrm{Pf}(4,W^\ast)_L)$ and $D(\mathrm{Pf}(4,W^\ast)_L)$, respectively.}. Moreover, given a generic linear subspace $L \subset \bigwedge^2W^\ast$ of dimension $7$, the linear section $X_L$ corresponds to a curve of genus $8$, and $Y_L$ is given by $\perf_\dg(\mathrm{Pf}(4,W^\ast)_L; \cF_L)$, where $\mathrm{Pf}(4,W^\ast)_L$ is a (singular) cubic $5$-fold; see \cite[\S10]{Lines}. Making use of Theorem \ref{thm:HPD}, we hence obtain equivalences of conjectures $C^+(X_L) \Leftrightarrow C^+_{\mathrm{nc}}(\perf_\dg(\mathrm{Pf}(4,W^\ast)_L; \cF_L))$ and $D(X_L) \Leftrightarrow D_{\mathrm{nc}}(\perf_\dg(\mathrm{Pf}(4,W^\ast)_L; \cF_L))$. Since Grothendieck's standard conjectures of type $C^+$ and $D$ are well-known in the case of curves, we hence conclude that the preceding noncommutative standard conjectures also hold.
\end{example}
As a second application, Theorem \ref{thm:HPD} shows us that whenever $X_L$, resp. $Y_L$, is of dimension $\leq 2$, the conjecture $C^+(Y_L)$, resp. $C^+(X_L)$, holds. Similarly, whenever $X_L$, resp. $Y_L$, is of dimension $\leq 4$, the conjecture $D(Y_L)$, resp. $D(X_L)$, holds. Here is an illustrative example; many more can be found at the survey \cite{ICM-Kuznetsov}:
\begin{example}[Determinantal duality]\label{ex:determinantal}
Let $U$ and $V$ be two $k$-vector spaces of dimensions $m$ and $n$, respectively, with $m\leq n$, $W$ the tensor product $U \otimes V$, and $0<r <m$ an integer. Following Bernardara-Bolognesi-Faenzi \cite[\S3]{BBF}, consider the determinantal variety $\mathcal{Z}^r_{m,n} \subset \bbP(W)$, resp. $\cW^r_{m,n}\subset \bbP(W^\ast)$, defined as the locus of those matrices $V \to U^\ast$, resp. $V^\ast \to U$, with rank at most $r$, resp. with corank at least $r$. For example, $\cZ^1_{m,n}$ are the classical Segre varieties. As explained in {\em loc. cit.}, $\mathcal{Z}^r_{m,n}$ and $\cW^r_{m,n}$ admit (Springer) resolutions of singularities $X:=\cX^r_{m,n}$ and $Y:=\cY^r_{m,n}$, respectively. Moreover, as proved in \cite[Thm.~3.5]{BBF}, $Y$ is the HP-dual of $X$. Making use of Theorem \ref{thm:HPD}, we hence conclude that $C^+(X_L)\Leftrightarrow C^+(Y_L)$ and $D(X_L) \Leftrightarrow D(Y_L)$ for every generic linear subspace $L \subset W$ of codimension $c$.
\end{example}
The linear section $X_L$ has dimension $r(m+n-r)-c-1$ and the linear section $Y_L$ dimension $r(m-n-r)+c-1$. Therefore, by combining the aforementioned work of Kleiman and Lieberman, with Example \ref{ex:determinantal}, we prove Grothendieck's standard conjectures in the case of linear sections of determinantal varieties:
\begin{theorem}[Linear sections of determinantal varieties]\label{thm:dimension}
Let $X_L$ and $Y_L$ be as in Example \ref{ex:determinantal}. Whenever $r(m+n-r)-c-1$ is $\leq 2$, resp. $\leq 4$, the conjecture $C^+(Y_L)$, resp. $D(Y_L)$, holds. Whenever $r(m-n-r)+c-1$ is $\leq2$, resp. $\leq 4$, the conjecture $C^+(X_L)$, resp. $D(X_L)$, holds.
\end{theorem}
%\begin{proof}
%As proved in \cite[\S3.4]{BBF}, all the assumptions of Theorem \ref{thm:HPD} are verified. Therefore, the proof follows from the combination of this latter theorem with the aforementioned work of Lieberman \cite{Lieberman}.
%\end{proof}
\begin{remark}[Dimension]
Note that Theorem \ref{thm:dimension} furnish us infinitely many examples of smooth projective $k$-schemes of arbitrary (high) dimension which satisfy Grothendieck's standard conjectures. For example, consider the case of square matrices, \ie $m=n$. Choose integers $n$, $r$, and $c$ (as above) such that $c<nr$ and $-r^2+c -1=2$. Under these choices, $Y_L$ has dimension $2$ and $X_L$ has dimension $2rn-(r^2 +c+1)>2$. Moreover, thanks to Theorem \ref{thm:dimension}, the conjectures $C^+(X_L)$ and $D(X_L)$ hold. Now, note that if we replce $n$ by $n'>n$ and keep $r$ and $c$, we obtain an higher dimensional $k$-scheme $X'_L$ of dimension $2rn'-(r^2+c+1)$ for which the conjectures $C^+(X'_L)$ and $D(X'_L)$ still hold.
\end{remark}

Theorem \ref{thm:main} allows us to easily extend Grothendieck's original conjectures from schemes $X$ to (smooth proper) stacks $\cX$ by setting $C^+(\cX):=C^+_{\mathrm{nc}}(\perf_\dg(\cX))$ and $D(\cX):=D_{\mathrm{nc}}(\perf_\dg(\cX))$. Making use of Theorem \ref{thm:HPD}, we now prove these extended conjectures in the case of bilinear divisors. Let $W$ be a $k$-vector space of dimension $d$, and $\cX$ the associated smooth proper Deligne-Mumford stack $(\bbP(W)\times \bbP(W))/\mu_2$ equipped with the map $\cX \to \bbP(S^2W), ([w_1],[w_2]) \mapsto [w_1\otimes w_2 + w_2 \otimes w_1]$. Given a generic linear subspace $L \subset S^2W^\ast$, the linear section $\cX_L$ corresponds to the intersection of the $\mathrm{dim}(L)$ bilinear divisors in $\cX$ parametrized by $L$.
\begin{theorem}[Intersections of bilinear divisors]\label{thm:last}
Assume that the dimension of $L$ is $\leq 3$, resp. $\leq 5$, and $d$ is odd or that the dimension of $L$ is $\leq 3$ and $d$ is even. In these cases, the conjecture $C^+(\cX_L)$, resp. $D(\cX_L)$, holds.
\end{theorem}
\begin{remark}
Voevodsky conjectured that the smash-nilpotence equivalence relation coincides with the numerical equivalence relation. %Since the homological equivalence relation sits between these two equivalence relations, Voevodsky's conjecture implies Grothendieck's conjecture. 
The corresponding analogues of Theorems \ref{thm:main}-\ref{thm:app1}, \ref{thm:HPD}, \ref{thm:app2}, and \ref{thm:dimension}, were established in \cite{Crelle}. %Note that since Voevodsky's conjecture is known only in the particular case where $\mathrm{dim}(X)\leq 2$, the aforementioned analogues are more restrictive.
\end{remark}

\section{Preliminaries}
%-------------------------------------------------------------------------------
%-------------------------------------------------------------------------------
\subsection{Dg categories}\label{sub:dg}
%-------------------------------------------------------------------------------
For a survey on dg categories consult Keller's ICM talk \cite{ICM-Keller}. Let $\cC(k)$ be the category of complexes of $k$-vector spaces. A {\em dg category $\cA$} is a category enriched over $\cC(k)$ and a {\em dg functor $F\colon \cA \to \cB$} is a functor enriched over $\cC(k)$. Let $\dgcat(k)$ be the category of (small) dg categories and dg functors. Recall from \cite[\S3.8]{ICM-Keller} that a {\em derived Morita equivalence} is a dg functor which induces an equivalence on derived categories. Following Kontsevich \cite{Miami,finMot,IAS}, a dg category $\cA$ is called {\em smooth} if it is perfect as a bimodule over itself and {\em proper} if $\sum_i \mathrm{dim}\, H^i\cA(x,y)< \infty$ for any pair of objects $(x,y)$. Examples include the dg categories of perfect complexes $\perf_\dg(X)$ associated to smooth proper $k$-schemes~$X$.

\subsection{Noncommutative motives}\label{sub:homological}
%-------------------------------------------------------------------------------
For a book on noncommutative motives consult \cite{book}. Recall from \cite[\S4.1]{book} the construction of the category of noncommutative Chow motives $\NChow(k)_\bbQ$. By construction, this rigid symmetric monoidal category comes equipped with a $\otimes$-functor $U(-)_\bbQ\colon \dgcat_{\mathrm{sp}}(k) \to \NChow(k)_\bbQ$ defined on smooth proper dg categories. Moreover, $\Hom_{\NChow(k)_\bbQ}(U(\cA)_\bbQ, U(\cB)_\bbQ)= K_0(\cA^\op \otimes \cB)_\bbQ$.
%As explained in \cite[\S 1.6.3]{book}, there is a natural bijection between $\Hom_{\Hmo(k)}(\cA,\cB)$ and the set of isomorphism classes of $\rep(\cA,\cB)$. Under this bijection, the composition law corresponds to the tensor product of dg bimodules. 
%The {\em additivization} of $\Hmo(k)$ is the additive category $\Hmo_0(k)$ with the same objects and with abelian groups of morphisms $\Hom_{\Hmo_0(k)}(\cA,\cB)$ given by the Grothendieck group $K_0\rep(\cA,\cB)$ of the triangulated category $\rep(\cA,\cB)$. The composition law is induced by the tensor product of bimodules. The {\em $\bbQ$-linearization} of $\Hmo_0(k)$ is the $\bbQ$-linear category $\Hmo_0(k)_\bbQ$ obtained by tensoring the morphisms of $\Hmo_0(k)$ with $\bbQ$. Note that we have the (composed) $\otimes$-functor
%\begin{eqnarray*}\label{eq:func3}
%U(-)_\bbQ\colon \dgcat(k) \too \Hmo_0(k)_\bbQ & \cA \mapsto \cA & (\cA\stackrel{F}{\to} \cB) \mapsto [{}_F\cB]_\bbQ\,.
%\end{eqnarray*}
%The category of {\em noncommutative Chow motives} $\NChow(k)_\bbQ$ is the idempotent completion of the full subcategory of $\Hmo_0(k)_\bbQ$ consisting of the objects $U(\cA)_\bbQ$ with $\cA$ smooth proper. This category is $F$-linear, additive, rigid symmetric monoidal, and idempotent complete; see \cite[\S4.1]{book}. Given dg categories $\cA$ and $\cB$, with $\cA$ smooth proper, we have $\rep(\cA, \cB) \simeq \cD_c(\cA^\op \otimes \cB)$. Hence, we obtain isomorphisms
%$$\Hom_{\NChow(k)_\bbQ}(U(\cA)_\bbQ, U(\cB)_\bbQ):=K_0(\rep(\cA,\cB))_\bbQ \simeq K_0(\cA^\op \otimes \cB)_\bbQ\,.$$
Recall from \cite[Thm.~9.2]{JEMS} that periodic cyclic homology gives rise to a $\bbQ$-linear $\otimes$-functor $HP^\pm\colon \NChow(k)_\bbQ \to \mathrm{Vect}_{\bbZ/2}(k)$ with values in finite dimensional $\bbZ/2$-graded $k$-vector spaces. The category of {\em noncommutative homological motives} $\NHom(k)_\bbQ$ is defined as the idempotent completion of the quotient $\NChow(k)_\bbQ/\mathrm{Ker}(HP^\pm)$. %By construction, $\NHom(k)_\bbQ$ is $\bbQ$-linear, additive, idempotent complete, and  rigid symmetric monoidal. 
Given a rigid symmetric monoidal category $\cC$, its $\cN$-ideal is defined as follows ($\mathrm{tr}(g\circ f)$ stands for the categorical trace of $g\circ f$):
$$\cN(a, b) := \{f \in \Hom_\cC(a, b) \,\,|\,\, \forall g \in \Hom_\cC(b, a)\,\, \mathrm{we}\,\,\mathrm{have}\,\, \mathrm{tr}(g\circ f) =0 \}\,.$$
Under these notations, the category of {\em noncommutative numerical motives} $\NNum(k)_\bbQ$ is defined as the idempotent completion of the quotient $\NChow(k)_\bbQ/\cN$. %By construction, $\NNum(k)_\bbQ$ is $\bbQ$-linear, additive, idempotent complete, and rigid symmetric monoidal.
%%-------------------------------------------------------------------------------
\subsection{Noncommutative standard conjecture of type $C^+$}\label{sub:C}
%%-------------------------------------------------------------------------------
Given a smooth proper dg category $\cA$, consider the K\"unneth projector $\pi_\cA^+\colon HP^\pm(\cA) \to HP^\pm(\cA)$. Following\footnote{In {\em loc. cit.} we used the notation $C_{\mathrm{nc}}(\cA)$ instead of $C^+_{\mathrm{nc}}(\cA)$.} \cite{JEMS}, the conjecture $C^+_{\mathrm{nc}}(\cA)$ asserts that $\pi^+_\cA$ is {\em algebraic}, \ie that there exists an endomorphism $\underline{\pi}^+_\cA \in \End_{\NChow(k)_\bbQ}(U(\cA)_\bbQ)$ such that $HP^\pm(\underline{\pi}^+_\cA)=\pi^+_\cA$.
%-------------------------------------------------------------------------------
\subsection{Noncommutative standard conjecture of type $D$}\label{sub:D}
%-------------------------------------------------------------------------------
Given a smooth proper dg category $\cA$, consider the $\bbQ$-vector spaces $K_0(\cA)_\bbQ/_{\!\sim\mathrm{hom}}$ and $K_0(\cA)_\bbQ/_{\!\sim\mathrm{num}}$ defined as $\Hom_{\NHom(k)_\bbQ}(U(k)_\bbQ, U(\cA)_\bbQ)$ and $\Hom_{\NNum(k)_\bbQ}(U(k)_\bbQ, U(\cA)_\bbQ)$. Following \cite{JEMS}, the conjecture $D(\cA)$ asserts that $K_0(\cA)_\bbQ/_{\!\sim\mathrm{hom}}= K_0(\cA)_\bbQ/_{\!\sim\mathrm{num}}$.
% are defined as follows:%equivalence relations on the $\bbQ$-linearization of the Grothendieck group $K_0(\cA)_\bbQ$ are defined as follows:
%\begin{eqnarray*}
%\Hom_{\NHom(k)_\bbQ}(U(k)_\bbQ, U(\cA)_\bbQ) && \Hom_{\NNum(k)_\bbQ}(U(k)_\bbQ, U(\cA)_\bbQ)\,.
%\end{eqnarray*}
%By construction, we have a surjective homomorphism $K_0(\cA)_\bbQ/_{\!\sim\mathrm{hom}} \twoheadrightarrow K_0(\cA)_\bbQ/_{\!\sim\mathrm{num}}$.
%The {\em standard conjecture $D_{\mathrm{nc}}(\cA)$} asserts that $K_0(\cA)_F/_{\!\sim\mathrm{hom}}= K_0(\cA)_F/_{\!\sim\mathrm{num}}$. 

%-------------------------------------------------------------------------------
\subsection{Orbit categories}\label{sub:orbit}
%-------------------------------------------------------------------------------
Let $(\cC, \otimes, {\bf 1})$ be an $\bbQ$-linear symmetric monoidal additive category and $\cO \in \cC$ a $\otimes$-invertible object. The {\em orbit category} $\cC/_{\!-\otimes \cO}$ has the same objects as $\cC$ and morphisms $\Hom_{\cC/_{\!-\otimes \cO}}(a,b):=\oplus_{n \in \bbZ} \Hom_\cC(a, b \otimes \cO^{\otimes n})$. Given objects $a, b, c$ and morphisms $\underline{f}=\{f_n\}_{n \in \bbZ}$ and $\underline{g}=\{g_n\}_{n \in \bbZ}$, the $i^{\mathrm{th}}$-component of $\underline{g}\circ \underline{f}$ is defined as $\sum_n (g_{i -n} \otimes \cO^{\otimes n})\circ f_n$. The canonical functor $\iota\colon \cC \to \cC/_{\!-\otimes \cO}$, given by $a \mapsto a$ and $f \mapsto \underline{f}=\{f_n\}_{n \in \bbZ}$, where $f_0=f$ and $f_n=0$ if $n\neq 0$, is endowed with an isomorphism $\iota \circ (-\otimes \cO) \Rightarrow \iota$ and is $2$-universal among all such functors. FInally, the category $\cC/_{\!-\otimes \cO}$ is $\bbQ$-linear, additive, and inherits from $\cC$ a symmetric monoidal structure making $\iota$ symmetric monoidal.
\section{Proof of Theorem \ref{thm:main}}
%-------------------------------------------------------------------------------
We start by proving the equivalence $C^+(X) \Leftrightarrow C^+_{\mathrm{nc}}(\perf_\dg(X))$. The implication $C^+(X) \Rightarrow C^+_{\mathrm{nc}}(\perf_\dg(X))$ was proved in \cite[Thm.~1.3]{JEMS}. Hence, we will prove solely the converse implication. Since $k$ is of characteristic zero, all the (classical) Weil cohomology theories $H^\ast$ are equivalent; see \cite[\S3.4.2]{Andre}. Therefore, in the proof we can (and will) make use solely of de Rham cohomology theory $H^\ast_{dR}$. 

Let us denote by $\Chow(k)_\bbQ$ the classical category of Chow motives. By construction, this rigid symmetric monoidal category comes equipped with a (contravariant) $\otimes$-functor $\mathfrak{h}(-)_\bbQ\colon \mathrm{SmProj}(k)^\op \to \Chow(k)_\bbQ$ defined on smooth projective $k$-schemes. As explained in \cite[Thm.~4.3]{book}, there exists a $\bbQ$-linear, fully faithful, $\otimes$-functor $\Phi$ making the following diagram commute
$$ 
\xymatrix{
\mathrm{SmProj}(k)^\op \ar[rr]^-{X\mapsto \perf_\dg(X)} \ar[d]_-{\mathfrak{h}(-)_\bbQ} && \dgcat_{\mathrm{sp}}(k) \ar[dd]^-{U(-)_\bbQ} \\
\Chow(k)_\bbQ \ar[d]_-\iota && \\
\Chow(k)_\bbQ/_{\!-\otimes \bbQ(1)} \ar[rr]_-{\Phi} && \NChow(k)_\bbQ\,,
}
$$
where $\Chow(k)_\bbQ/_{\!-\otimes \bbQ(1)}$ stands for the orbit category with respect to the Tate motive $\bbQ(1)$. Consider the following composition:
 \begin{equation}\label{eq:composition}
 \Chow(k)_\bbQ \stackrel{\iota}{\too} \Chow(k)_\bbQ/_{\!-\otimes \bbQ(1)} \stackrel{\Phi}{\too} \NChow(k)_\bbQ \stackrel{HP^\pm}{\too} \mathrm{Vect}_{\bbZ/2}(k)\,.
 \end{equation}
Consider also the category $\mathrm{Vect}_\bbZ(k)$ of finite dimensional $\bbZ$-graded $k$-vector spaces and the associated ``$2$-perioditization'' $\otimes$-functor:
\begin{eqnarray}\label{eq:perioditization}
\mathrm{Vect}_\bbZ(k) \too \mathrm{Vect}_{\bbZ/2}(k) && \{V_i\}_{i \in \bbZ} \mapsto (\oplus_{i\,\mathrm{even}} V_i, \oplus_{i\,\mathrm{odd}} V_i)\,. 
\end{eqnarray} 
Given any smooth (projective) $k$-scheme $X$, the Hochschild-Kostant-Rosenberg theorem identifies $HP^+(\perf_\dg(X))$ and $HP^-(\perf_\dg(X))$ with $\oplus_{i\,\mathrm{even}} H^i_{dR}(X)$ and $\oplus_{i\,\mathrm{odd}} H^i_{dR}(X)$, respectively.  Consequently, \eqref{eq:composition} reduces to the composition of 
 \begin{eqnarray}\label{eq:deRham}
 H^\ast_{dR}\colon \Chow(k)_\bbQ \too \mathrm{Vect}_{\bbZ}(k) && \mathfrak{h}(X)_\bbQ \mapsto  \oplus_{i=0}^{2\mathrm{dim}(X)} H^i_{dR}(X)
 \end{eqnarray}
with the above functor \eqref{eq:perioditization}. Assume now conjecture $C^+_{\mathrm{nc}}(\perf_\dg(X))$. Since $\Phi$ is an equivalence of categories, there exists then an endomorphism $\underline{\pi}^+=\{\pi^+_n\}_{n \in \bbZ}$ of $\iota(\mathfrak{h}(X)_\bbQ)$ in the orbit category $\Chow(k)_\bbQ/_{\!-\otimes \bbQ(1)}$ whose image under the composed functor $HP^\pm \circ \Phi$ agrees with the homomorphism of $\bbZ/2$-graded $k$-vector spaces:
$$ \pi^+\colon (\oplus_{i\,\mathrm{even}} H^i_{dR}(X), \oplus_{i\,\mathrm{odd}} H^i_{dR}(X)) \stackrel{(\id,0)}{\too} (\oplus_{i\,\mathrm{even}} H^i_{dR}(X), \oplus_{i\,\mathrm{odd}} H^i_{dR}(X))\,.$$
Note that $\pi^+$ is the image of the even K\"unneth projector $\pi_X^+:=\sum_i \pi_X^{2i}$ under~\eqref{eq:perioditization}. Note also that $H^\ast_{dR}(\pi_n^+)$ is an homomorphism of degree $-2n$. The preceding considerations, combined with the construction of the orbit category $\Chow(k)_\bbQ/_{\!-\otimes \bbQ(1)}$, allows us then to conclude that $H^\ast_{dR}(\pi^+_0)=\pi_X^+$ and $H^\ast_{dR}(\pi^+_n)=0$ if $n\neq 0$. Consequently, the even K\"unneth projector $\pi_X^+$ is algebraic and conjecture $C^+(X)$ holds.

Let us now prove the equivalence $D(X) \Leftrightarrow D_{\mathrm{nc}}(\perf_\dg(X))$. The implication $D(X) \Rightarrow D_{\mathrm{nc}}(\perf_\dg(X))$ was proved in \cite[Thm.~1.5]{JEMS}. Hence, we will prove solely the converse implication. Recall from \cite[page~645]{JEMS}, the construction of the following commutative square:
\begin{equation}\label{eq:square1}
\xymatrix{
Z^\ast(X)_\bbQ/_{\!\sim \mathrm{hom}} \ar@{->>}[r] \ar@{->>}[d] & K_0(\perf_\dg(X))_\bbQ/_{\!\sim \mathrm{hom}} \ar@{->>}[d] \\
Z^\ast(X)_\bbQ/_{\!\sim \mathrm{num}} \ar[r]_-\simeq & K_0(\perf_\dg(X))_\bbQ/_{\!\sim \mathrm{num}}\,.
}
\end{equation}
Note that the conjecture $D(X)$, resp. $D_{\mathrm{nc}}(\perf_\dg(X))$, is equivalent to the injectivity of the vertical left-hand side, resp. right-hand side, homomorphism. Assume now conjecture $D_{\mathrm{nc}}(\perf_\dg(X))$, \ie that the vertical right-hand side homomorphism in \eqref{eq:square1} is injective. By construction, the vertical left-hand side homomorphism in \eqref{eq:square1} is a diagonal (matrix) homomorphism:
$$ Z^\ast(X)_\bbQ/_{\!\sim \mathrm{hom}} = \oplus_{n=0}^{\mathrm{dim}(X)} Z^n(X)_\bbQ/_{\!\sim \mathrm{hom}} \to \oplus_{n=0}^{\mathrm{dim}(X)} Z^n(X)_\bbQ/_{\!\sim \mathrm{num}} = Z^\ast(X)_\bbQ/_{\!\sim \mathrm{num}}\,.$$
Therefore, in order to prove conjecture $D(X)$ it suffices to show that the following homomorphisms are injective:
\begin{eqnarray}\label{eq:homomorphisms}
Z^n(X)_\bbQ/_{\!\sim \mathrm{hom}} \too K_0(\perf_\dg(X))_\bbQ/_{\!\sim \mathrm{hom}} && 0 \leq n \leq \mathrm{dim}(X)\,.
\end{eqnarray}
As explained above, we have $\eqref{eq:composition}=\eqref{eq:perioditization}\circ \eqref{eq:deRham}$. This implies that the classical category of homological motives $\Hom(k)_\bbQ$ agrees with the idempotent completion of the quotient $\Chow(k)_\bbQ/\mathrm{Ker}(\eqref{eq:composition})$. Moreover, the induced functor $\Hom(k)_\bbQ  \to \NHom(k)_\bbQ$ is faithful. Under the following identifications
$$ \Hom_{\Hom(k)_\bbQ}(\mathfrak{h}(\mathrm{Spec}(k))_\bbQ, \mathfrak{h}(X)_\bbQ(n)) \simeq Z^n(X)_\bbQ/_{\!\sim \mathrm{hom}}$$
$$ \Hom_{\NHom(k)_\bbQ}(U(k)_\bbQ, U(\perf_\dg(X))_\bbQ \simeq K_0(\perf_\dg(X))_\bbQ/_{\!\sim \mathrm{hom}}\,,$$
where $\mathfrak{h}(X)_\bbQ(n)$ stands for $\mathfrak{h}(X)_\bbQ \otimes \bbQ(1)^{\otimes n}$, the homomorphisms \eqref{eq:homomorphisms} correspond to the homomorphisms induced by the functor $\Hom(k)_\bbQ  \to \NHom(k)_\bbQ$:
$$\Hom_{\Hom(k)_\bbQ}(\mathfrak{h}(\mathrm{Spec}(k))_\bbQ, \mathfrak{h}(X)_\bbQ(n)) \too \Hom_{\NHom(k)_\bbQ}(U(k)_\bbQ, U(\perf_\dg(X))_\bbQ)\,.$$
Since these latter homomorphisms are injective, we hence conclude that the conjecture $D(X)$ holds. This finishes the proof.
%-------------------------------------------------------------------------------
\section{Proof of Theorem \ref{thm:app1}}
%-------------------------------------------------------------------------------
As proved in \cite[Thm.~4.2]{Kuznetsov-quadrics}, the category $\perf(Q)$ admits a semi-orthogonal decomposition $\langle \perf(S;\cF), \perf(S)_1, \ldots, \perf(S)_d\rangle$, where $\cF$ stands for the sheaf of even Clifford algebras associated to $q$ and $\perf(S)_i:=q^\ast \perf(S) \otimes \cO_{Q/S}(i)$. Note that $\perf(S)_i\simeq \perf(S)$. As explained in \cite[\S2.4.1]{book}, this semi-orthogonal decomposition gives rise to a direct sum decomposition in the category $\NChow(k)_\bbQ$:
$$ U(\perf_\dg(Q))_\bbQ \simeq U(\perf_\dg(S;\cF))_\bbQ \oplus U(\perf_\dg(S))_\bbQ \oplus \cdots \oplus U(\perf_\dg(S))_\bbQ\,.$$
Making use of the definition of the noncommutative standard conjectures of type $C^+$ and $D$, we hence obtain the following equivalence of conjectures:
\begin{equation}\label{eq:conj0}
C^+_{\mathrm{nc}}(\perf_\dg(Q)) \Leftrightarrow C^+_{\mathrm{nc}}(\perf_\dg(S;\cF)) + C^+_{\mathrm{nc}}(\perf_\dg(S))
\end{equation}
\begin{equation}\label{eq:conj1}
D_{\mathrm{nc}}(\perf_\dg(Q)) \Leftrightarrow D_{\mathrm{nc}}(\perf_\dg(S;\cF)) + D_{\mathrm{nc}}(\perf_\dg(S))\,.
\end{equation}
Since $d$ is even and the discriminant division of $q$ is smooth, the category $\perf(S;\cF)$ is equivalent (via a Fourier-Mukai functor) to $\perf(\widetilde{S};\widetilde{\cF})$ where $\widetilde{S}$ is the discriminant double cover of $S$ and $\widetilde{\cF}$ a sheaf of Azumaya algebras over $\widetilde{S}$; see \cite[Prop.~3.13]{Kuznetsov-quadrics}. This implies that $\perf_\dg(S;\cF)$ is derived Morita equivalent to $\perf_\dg(\widetilde{S};\widetilde{\cF})$. Using the fact that $U(-)_\bbQ$ inverts Morita equivalences and the isomorphism between $U(\perf_\dg(\widetilde{S};\widetilde{\cF}))_\bbQ$ and $U(\perf_\dg(\widetilde{S}))_\bbQ$ established in \cite[Thm.~2.1]{Azumaya}, we hence conclude that the right-hand side of \eqref{eq:conj0}, resp. \eqref{eq:conj1}, reduces to $C^+_{\mathrm{nc}}(\perf_\dg(\widetilde{S}))+ C^+_{\mathrm{nc}}(\perf_\dg(S))$, resp. to $D_{\mathrm{nc}}(\perf_\dg(\widetilde{S}))+ D_{\mathrm{nc}}(\perf_\dg(S))$. Finally, since the dimension of $S$ is $\leq 2$, resp. $\leq 4$, the aforementioned work of Kleiman and Lieberman, combined with Theorem \ref{thm:main}, implies conjecture $C^+(Q)$, resp. $D(Q)$.
%-------------------------------------------------------------------------------
\section{Proof of Theorem \ref{thm:HPD}}
%-------------------------------------------------------------------------------
The proof is similar for the noncommutative standard conjecture of type $C^+$ and type $D$. Therefore, we will prove solely the first case.

By definition of the Lefschetz decomposition $\langle \bbA_0, \bbA_1(1), \ldots, \bbA_{i-1}(i-1)\rangle$, we have a chain of admissible triangulated subcategories $\bbA_{i-1}\subseteq \cdots \subseteq \bbA_1\subseteq \bbA_0$ and $\bbA_r(r):=\bbA_r \otimes \cO_X(r)$. Note that $\bbA_r(r)\simeq \bbA_r$. Let $\mathfrak{a}_r$ be the right orthogonal complement to $\bbA_{r+1}$ in $\bbA_r$; these are called the {\em primitive subcategories} in \cite[\S4]{KuznetsovHPD}. Note that we have semi-orthogonal decompositions:
\begin{eqnarray}\label{eq:decomp1}
\bbA_r = \langle \mathfrak{a}_r, \mathfrak{a}_{r+1}, \ldots, \mathfrak{a}_{i-1} \rangle && 0\leq r \leq i-1\,.
\end{eqnarray}
As proved in \cite[Thm.~6.3]{KuznetsovHPD}, the category $\perf(Y)$ admits a HP-dual Lefschetz decomposition $\perf(Y)=\langle \bbB_{j-1}(1-j), \bbB_{j-2}(2-j), \ldots, \bbB_0\rangle$ with respect to $\cO_Y(1)$. As above, we have a chain of admissible subcategories $\bbB_{j-1} \subseteq \bbB_{j-2} \subseteq \cdots \subseteq \bbB_0$. Moreover, the primitive subcategories coincide (via a Fourier-Mukai functor) with those of $\perf(X)$ and we have semi-orthogonal decompositions:
\begin{eqnarray}\label{eq:decomp2}
\bbB_r=\langle \mathfrak{a}_0, \mathfrak{a}_1, \ldots, \mathfrak{a}_{\mathrm{dim}(V)-r-2}\rangle && 0 \leq r \leq j-1\,.
\end{eqnarray}
Furthermore, the assumptions $\mathrm{dim}(X_L)=\mathrm{dim}(X)-\mathrm{dim}(L)$ and $\mathrm{dim}(Y_L)=\mathrm{dim}(Y) - \mathrm{dim}(L^\perp)$ imply the existence of semi-orthogonal decompositions
\begin{equation}\label{eq:semi-1}
\perf(X_L)=\langle \bbC_L, \bbA_{\mathrm{dim}(V)}(1), \ldots, \bbA_{i-1}(i-\mathrm{dim}(V))\rangle 
\end{equation}
\begin{equation}\label{eq:semi-2}
\perf(Y_L)=\langle \bbB_{j-1}(\mathrm{dim}(L^\perp)-j), \ldots, \bbB_{\mathrm{dim}(L^\perp)}(-1), \bbC_L \rangle\,,
\end{equation}
where $\bbC_L$ is a common triangulated category. Let us denote by $\bbC_L^\dg$, $\bbA_r^\dg$ and $\mathfrak{a}_r^\dg$ the dg enhancement of $\bbC_L$, $\bbA_r$ and $\mathfrak{a}_r$ induced from $\perf_\dg(X_L)$. Similarly, let us denote by $\bbC_L^{\dg'}$ and $\bbB_r^{\dg}$ the dg enhancement of $\bbC_L$ and $\bbB_r$ induced from $\perf_\dg(Y_L)$. Note that since $X_L$ and $Y_L$ are smooth projective $k$-schemes, all the preceding dg categories are smooth proper; see \cite[Lem.~2.1]{Crelle}. As explained in \cite[\S2.4.1]{book}, the above semi-orthogonal decompositions \eqref{eq:semi-1}-\eqref{eq:semi-2} give rise to the following direct sums decompositions in the category $\NChow(k)_\bbQ$:
$$ U(\perf_\dg(X_L))_\bbQ \simeq U(\bbC^\dg_L)_\bbQ \oplus U(\bbA^\dg_{\mathrm{dim}(V)})_\bbQ \oplus \cdots \oplus U(\bbA^\dg_{i-1})_\bbQ$$
$$ U(\perf_\dg(Y_L))_\bbQ \simeq U(\bbB_{j-1}^\dg)_\bbQ \oplus \cdots \oplus U(\bbB^\dg_{\mathrm{dim}(L^\perp)})_\bbQ \oplus U(\bbC_L^{\dg'})_\bbQ\,.$$
Making use of the definition of the noncommutative standard conjecture of type $C^+$, we hence obtain the following equivalences of conjectures:
\begin{equation}\label{eq:D1}
C^+_{\mathrm{nc}}(\perf_\dg(X_L))\Leftrightarrow C^+_{\mathrm{nc}}(\bbC_L^\dg)+ C^+_{\mathrm{nc}}(\bbA^\dg_{\mathrm{dim}(V)}) + \cdots + C^+_{\mathrm{nc}}(\bbA^\dg_{i-1})
\end{equation}
\begin{equation}\label{eq:D2}
C^+_{\mathrm{nc}}(\perf_\dg(Y_L))\simeq C^+_{\mathrm{nc}}(\bbB^\dg_{j-1}) + \cdots + C^+_{\mathrm{nc}}(\bbB^\dg_{\mathrm{dim}(L^\perp)}) + C^+_{\mathrm{nc}}(\bbC^{\dg'}_L)\,.
\end{equation}
On the one hand, since the conjecture $C^+_{\mathrm{nc}}(\bbA_0^\dg)$ holds, we conclude from the semi-orthogonal decompositions \eqref{eq:decomp1}-\eqref{eq:decomp2} that the conjectures $C^+_{\mathrm{nc}}(\bbA^\dg_r)$ and $C^+_{\mathrm{nc}}(\bbB^\dg_r)$ hold for every $r$. This implies that the right-hand side of \eqref{eq:D1}, resp. \eqref{eq:D2}, reduces to $C^+_{\mathrm{nc}}(\bbC^\dg_L)$, resp. $C^+_{\mathrm{nc}}(\bbC_L^{\dg'})$. On the other hand, since the composed functor $\perf(X_L) \to \bbC_L \to \perf(Y_L)$ is of Fourier-Mukai type, the dg categories $\bbC_L^\dg$ and $\bbC_L^{\dg'}$ are derived Morita equivalent. Using the fact that the functor $U(-)_\bbQ$ inverts derived Morita equivalences, this implies that $C^+_{\mathrm{nc}}(\bbC^\dg_L)\Leftrightarrow C^+_{\mathrm{nc}}(\bbC_L^{\dg'})$. Finally, since $X_L$ and $Y_L$ are smooth (projective) $k$-schemes, the proof follows from Theorem \ref{thm:main}.
\section{Proof of Theorem \ref{thm:app2}}
%-------------------------------------------------------------------------------
Similarly to the proof of Theorem \ref{thm:app1}, since the fibration $q\colon Q\to \bbP(S^2W^\ast)$ is of relative dimension $d-2$, $d$ is even, and the discriminant divisor of $q$ is smooth, the conjecture $C^+_{\mathrm{nc}}(\perf_\dg(\bbP(L);\cF_L))$, resp. $D_{\mathrm{nc}}(\perf_\dg(\bbP(L);\cF_L))$, reduces to conjecture $C^+_{\mathrm{nc}}(\perf_\dg(\widetilde{\bbP(L)}))$, resp.  $D_{\mathrm{nc}}(\perf_\dg(\widetilde{\bbP(L)}))$. The proof follows then from the assumption that the dimension of $L$ is $\leq 3$, resp. $\leq 5$, from the aforementioned work of Kleiman and Lieberman, and from Theorem \ref{thm:main}.

%-------------------------------------------------------------------------------
\section{Proof of Theorem \ref{thm:last}}
%-------------------------------------------------------------------------------
As proved by Rennemo in \cite[Thm.~1.1 and Props.~1.2-1.3]{Rennemo}, in the cases of Theorem \ref{thm:last} the linear section $Y_L$ of the ``HP-dual'' of $\cX$ is given by a certain smooth double cover of $\bbP(L)$. Consequently, since the dimension of $Y_L$ is $\leq 2$, resp. $\leq 4$, the proof follows from the combination of the aforementioned work of Kleiman and Lieberman, with Theorem~\ref{thm:HPD}.

\medbreak\noindent\textbf{Acknowledgments:} I thank Bruno Kahn for reminding me that Grothendieck's standard conjecture of Lefschetz type is stable under hyperplane sections.

\end{document}

\end{proof}